# ABOUT CONVERGENCE OF PROJECTION METHODS FOR SOLUTION OF SOME FREDHOLM INTEGRAL EQUATIONS OF THE FIRST KIND

Olexandr Polishchuk

Pidstryhach Institute for Applied Problems of Mechanics and Mathematics,
National Academy of Sciences of Ukraine, Lviv, Ukraine

od_polishchuk@ukr.net

**Abstract.** This article is dedicated to research of approximation properties of B-splines and Lagrangian finite elements in Hilbert spaces of functions defined on surfaces in three-dimensional space. Hereinafter the conditions are determined for convergence of Galerkin and collocation methods for solving Fredholm integral equation of the first kind for simple layer potential that is equivalent to Dirichlet problem for Laplace equation in $R^3$. Estimation is determined for the error of approximate solution of this problem obtained using potential theory methods.

## 1. Introduction

Many physical processes (e.g. diffusion, heat flux, electrostatic field, perfect fluid flow, elastic motion of solid bodies, groundwater flow, etc.) are modeled using boundary value problems for Laplace equation [1–3]. The powerful tools for solving such problems are potential theory methods, especially in the case of tired boundary surface or complex shape surface [4–6]. These methods are a convenient for calculating desired solution in small domains [7, 8]. In number of cases, application of potential theory methods requires solving Fredholm integral equation of the first kind. In particular, one of the cases is solving Dirichlet problem in the space of functions with normal derivative jump on crossing boundary surface using simple layer potential [9, 10]. When solving Neumann problem in the space of functions with jump on crossing boundary surface using double layer potential, we also proceed to integral equation of the first kind [11, 12]. The need to solve integral equations of the first kind also arises when the sum of simple and double layer potentials is used to solve the double-sided Dirichlet or Neumann problem [13] or double-sided Dirichlet-Neumann problem [14] in the space of functions that, same as their normal derivatives, have jump on crossing boundary surface. Many systems of integral equations for the simple and double layer potentials that are equivalent to mixed boundary value problems for Laplace equation, also contain integral equations of the first kind [15, 16].

In general, researches of projection methods convergence mainly focus on solving integral equations of the second kind [4, 6, 17]. Defining well-posed solvability conditions for integral equations of the first kind that are equivalent to boundary value problems for Laplace's equation in Hilbert spaces [18–20] allows us to use projection methods for numerical solution of such equations, thus avoiding resource-consuming regularization procedures [21–23]. For detailed review of numerical methods for solving integral equations, please see [2–4, 6]. In [24, 25], convergence conditions are defined for the series of projection methods for solving Fredholm integral equation of the first kind for simple layer potential that is equivalent to three-dimensional Dirichlet problem for Laplace equation while approximating desired potential density with complete systems of orthonormal functions. However, if boundary surface has a complex shape usage of such approximations poses considerable difficulties for practical implementation of numerical methods [7]. In this case, finite elements of different types should be used for approximation of desired potential densities [26, 27]. Derived approximations, among other things, allow us to create effective algorithms for singularities removal in kernels and desired integral equation densities [28].

The purpose of the paper is to define convergence conditions of projection methods for approximate solution of Fredholm integral equations of the first kind by the example of integral equation for the simple layer potential that is equivalent to Dirichlet problem for Laplace equation using approximation of desired potential density with systems of finite elements of different types and orders.

## 2. Approximations of Hilbert spaces and the basic convergence theorem

Consider the operator equation
$$Au = f, \ u \in U, \ f \in F, \qquad (1)$$
where $U$ and $F$ are the Hilbert spaces, $A \in L(U, F)$. To solve equation (1) we apply the projection method
$$Q_N A P_N u = Q_N f. \qquad (2)$$

In formula (2) $P_N$ and $Q_N$ are projection operators from $U$ and $F$ onto closed finite dimensional subspaces $U_N \in U$ and $F_N \in F$ accordingly. Define operators $P_N$ and $Q_N$ in the following way. Denote by $r_N$ the restriction operator from the space $U$ to the finite dimensional subspace $V_N \subset R^N$ and introduce in $V_N$ the extension operator $p_N$ as an isomorphism from $V_N$ onto subspace $U_N \in U$. The norm in $V_N$ is determined by the relation $\|\mathbf{u}_N\|_{V_N} = \|p^N \mathbf{u}_N\|_U$, $\mathbf{u}_N \in V_N$. Then $P_N = p_N r_N$ and we can determine the triple $(V_N, p_N, r_N)$ as approximation of the space $U$. Such approximation are called convergent if
$$\lim_{N \to \infty} \|u - p_N r_N u\|_U = 0.$$

Denote by $s_N$ the restriction operator from the space $F$ to the finite dimensional subspace $\Phi_N \subset R^N$. The extension operator $q_N$ from $\Phi_N$ onto subspace $F_N = AU_N \subset F$ introduce by the formula
$$q_N \mathbf{f}_N = A p_N \mathbf{u}_N, \ \mathbf{f}_N \in \Phi_N.$$
The norm in $\Phi_N$ is determined by the relation $\|\mathbf{f}_N\|_{\Phi_N} = \|q^N \mathbf{f}_N\|_F$, $\mathbf{f}_N \in \Phi_N$. Then $Q_N = q_N s_N$ and we can determine the triple $(\Phi_N, q_N, s_N)$ as approximation of the space $F$. Thus, the solution of problem (2) is reduced to solution of the system of linear algebraic equations
$$\mathbf{A}_N \mathbf{u}_N = \mathbf{f}_N, \ \mathbf{A}_N = s_N A p_N, \ \mathbf{A}_N \in L(V_N, \Phi_N), \ \mathbf{f}_N = s_N f. \qquad (3)$$

Operators $\mathbf{A}_N$ are called stable if there is independent of $N$ constant $\mu > 0$ such that for arbitrary $\mathbf{u}_N \in V_N$ is performed inequality
$$\mu \|\mathbf{u}_N\|_{V_N} \leq \|\mathbf{A}_N \mathbf{u}_N\|_{F_N}. \qquad (4)$$

Let us the pairs of operators $(r_N, p_N)$ and $(s_N, q_N)$ are selected. Assume as an approximate solution of equation (1) the function $p_N \mathbf{u}_N \in U_N$ where $\mathbf{u}_N$ is the solution of problem (3). Then we have the next basic theorem of convergence [29].

Theorem 1. Let us the operator $A$ is an isomorphism from $U$ into $F$. Then the sequence $p_N \mathbf{u}_N$ converges to solution $u$ of equation (1) if and only if the approximations $(V_N, p_N, r_N)$ of the space $U$ are convergent and operators $\mathbf{A}_N$ are stable. In addition, error estimation of the approximate solution is given by the ratio
$$\|u - p_N \mathbf{u}_N\|_U \leq (1 + \|A\|/\mu) \|u - p_N r_N u\|_U.$$

The choice of triples $(V_N, p_N, r_N)$ and $(\Phi_N, q_N, s_N)$ defines one or another projection method of solving the equation (1) (Galerkin, smallest squares, smallest mismatch, collocation, etc.).

**3. B-splines approximations**

Let us $S = [0,a] \times [0,b] \subset R^2$. Construct in the domain $S$ a rectangular grid $S_h$ with the steps $h_1 = a/n$ and $h_2 = b/k$, $n,k = 1,2,...$ Introduce in the space $H^m(S)$, $m = 0,1,2,...$, the system of B-splines of $m$-th degree

$$\{B_{ij}(\xi)\}_{i=-m}^{n-1} {}_{j=-m}^{k-1}, \; n,k \geq m+1. \tag{5}$$

Denote by $U_B^N$ the linear shell of system (5). Select restriction operator $r_B^N : H^m(S) \to V_B^N \subset R^N$ and extension operator $p_B^N : V_B^N \to U_B^N \subset H^m(S)$ in the form

$$(r_B^N v)_{i,j} \equiv v_N^{(i,j)} = \int_S B_{ij}(\xi) v(\xi) dS_\xi, \; i = \overline{-m(1)(n-1)}, \; j = \overline{-m(1)(k-1)}, \tag{6}$$

$$r_B^N v = \mathbf{v}_B^N \in V_B^N, \; v \in H^m(S),$$

$$p_B^N \mathbf{v}_B^N = \sum_{i=-m}^{n-1} \sum_{j=-m}^{k-1} v_N^{(i,j)} B_{ij}(\xi), \; N = (n+m)(k+m). \tag{7}$$

The next result is in order [29].

Lemma 1. Approximations $(V_B^N, p_B^N, r_B^N)$ of the space $H^m(S)$ are convergent and for arbitrary $v \in H^m(S)$ are performed the estimates

$$\|v - p_B^N r_B^N v\|_{H^t(S)} \leq C h^{\sigma-t} \|v\|_{H^\sigma(S)}, \; 0 \leq t \leq \sigma \leq m+1, \; t \leq m, \tag{8}$$

where constant $C > 0$ does not depend on $v$.

Denote by $G$ a bounded open domain in $R^3$ with boundary $\Gamma$. Suppose that exists $M$ open balls

$$B_l \subset R^3, \; \Gamma \subset \bigcup_{l=1}^M B_l, \; B_l \cap \Gamma = \Gamma_l \neq 0, \; l = \overline{1,M},$$

such that for each ball $B_l$ there is defined on $\overline{B}_l$ $m$ times differentiated real vector-function $f^{(l)}(x) = (f_1^{(l)}, f_2^{(l)}, f_3^{(l)})$ such that $y_l = f^{(l)}(x)$ carries a mutually unambiguous mapping of the ball $B_l$ onto some open bounded set in $R^3$ where $\Gamma_l$ is mapped on an open set $S_l \subset R^2$. In addition, the Jacobian $J_l = \dfrac{\partial(f_1^{(l)}, f_2^{(l)}, f_3^{(l)})}{\partial(x_1, x_2, x_3)}$ is positive and continuous if $x \in \overline{B}_l$, $l = \overline{1,M}$. Then the surface $\Gamma$ is called m-smooth surface in $R^3$ [30].

We associate to the partition $\Gamma = \bigcup_{l=1}^M \Gamma_l$ the partition of one $\{\psi_l(x)\}_{l=1}^M$, $x \in \Gamma$, with the following properties:

$$\psi_l(x) \in C^\infty(\Gamma_l), \; \text{supp}\{\psi_l\} \subset \Gamma_l, \; 0 \leq \psi_l(x) \leq 1, \; l = \overline{1,M}, \; \sum_{l=1}^M \psi_l(x) = 1,$$

and exists $m$ times continuously differentiated mapping

$$\tau_l : \Gamma_l \to S_l = [0,a_l] \times [0,b_l], \; \tau_l^{-1} : S_l \to \Gamma_l, \; l = \overline{1,M}.$$

Then for arbitrary function $u(x)$ defined on $\Gamma$ we can put into a mutually unambiguous correspondance the set of defined in $R^2$ functions

$$\{v_l(\xi)\}_{l=1}^{M},\ v_l(\xi)=u(\tau_l^{-1}(\xi)),\ \xi\in S_l,$$

which have the compact support on $S_l$ and

$$u_l(x)=\psi_l(x)u(x),\ x\in\Gamma,\ u(x)\in H^m(\Gamma),\ \text{if}\ v_l(\xi)\in H^m(S_l),\ l=\overline{1,M},$$

and

$$\|u\|_{H^m(\Gamma)}=\sum_{l=1}^{M}\|v_l\|_{H^m(S_l)}.$$

Construct in each domain $S_l$ a rectangular grid $S_l^h$ with the steps $h_1^{(l)}=a_l/n_l$ and $h_2^{(l)}=b_l/k_l$, and define in each grid domain $S_l^h$ a system of functions

$$\{B_{ij}^{(l)}(\xi)\}_{i=-m,\ j=-m}^{n_l-1,\ k_l-1},\ n_l,k_l\geq m+1,\ l=\overline{1,M}.$$

Assign to them the grid $\Gamma_h=\bigcup_{l=1}^{M}\tau_l^{-1}(S_l^h)$ on the surface $\Gamma$ and system of functions

$$\{\widetilde{B}_k(x)\}_{k=1}^{N_B}=\bigcup_{l=1}^{M}\{B_{ij}^{(l)}(\tau_l(x))\}_{i=-m,\ j=-m}^{n_l-1,\ k_l-1},\ N_B=\sum_{l=1}^{M}(n_l+m)(k_l+m). \tag{9}$$

Denote $\Gamma_{l_1,l_2}=\Gamma_{l_1}\cap\Gamma_{l_2}$, $l_1,l_2=\overline{1,M}$, and suppose that the grid on surface $\Gamma$ satisfies condition

$$\mathrm{supp}\{B_{ij}^{(l)}(\tau_l(x))\}\not\subset\Gamma_{l,k},\ i=\overline{-m(1)(n_l-1)},\ j=\overline{-m(1)(k_l-1)},\ k\neq l,\ k,l=\overline{1,M}.$$

Since $\mathrm{supp}\{\widetilde{B}_{i_1}(x)\}\neq\mathrm{supp}\{\widetilde{B}_{i_2}(x)\}$ for $i_1\neq i_2$, the functions of system $\{\widetilde{B}_k(x)\}_{k=1}^{N_B}$ are linearly independent.

Denote by $\widetilde{r}_B^{N_B}$ the restriction operator from $H^m(\Gamma)$ onto finite dimensional space $V_B^{N_B}$, and $\widetilde{r}_B^{N_l}$ is its restriction to $H^m(\Gamma_l)$, that is,

$$\widetilde{r}_B^{N_B}=\{\widetilde{r}_B^{N_l}\}_{l=1}^{M},\ \widetilde{r}_B^{N_l}u_l(x)=r_B^{N_l}v_l(\xi),\ N_l=(n_l+m)(k_l+m), \tag{10}$$

where $r_B^{N_l}$ is the similar to (6) restriction operator from $H^m(S_l)$ onto finite dimensional space $V_B^{N_l}$, $l=\overline{1,M}$, and $V_B^{N_B}=V_B^{N_1}\times V_B^{N_2}\times...\times V_B^{N_M}$.

The extension operator $\widetilde{p}_B^{N_B}$ from $V_B^{N_B}$ onto $U_B^{N_B}\subset H^m(\Gamma)$ is introduced by the formula

$$(\widetilde{p}_B^{N_B}\mathbf{u}_{N_B})(x)=\sum_{i=1}^{N_B}u_i\widetilde{B}_i(x),\ \mathbf{u}_{N_B}\in V_B^{N_B}. \tag{11}$$

From lemma 1 follows that

$$\lim_{N_B\to\infty}\|u-\widetilde{p}_B^{N_B}\widetilde{r}_B^{N_B}u\|_{H^m(\Gamma)}=\sum_{l=1}^{M}\lim_{N_l\to\infty}\|v_l-p_B^{N_l}r_B^{N_l}v_l\|_{H^m(S_l)}=0,$$

i.e. approximations $(V_B^{N_B},\widetilde{p}_B^{N_B},\widetilde{r}_B^{N_B})$ of the space $H^m(\Gamma)$ are convergent. Further, from estimates (8) we obtain

$$\|u-\widetilde{p}_B^{N_B}\widetilde{r}_B^{N_B}u\|_{H^t(\Gamma)}^2=\sum_{l=1}^{M}\|v_l-p_B^{N_l}r_B^{N_l}v_l\|_{H^t(S_l)}^2\leq C^2h^{2(\sigma-t)}\sum_{l=1}^{M}\|v_l\|_{H^\sigma(S_l)}^2=C^2h^{2(\sigma-t)}\|u\|_{H^\sigma(\Gamma)}^2,$$

$$0\leq t\leq\sigma\leq m+1,\ t\leq m,$$

where $p_B^{N_l}$ is the similar to (7) extension operator from $V_B^{N_l}$ onto the spaces $H^m(S_l)$, constant $C>0$ does not depend on $u$, and $h=\max_{1\leq l\leq M}\{h_1^{(l)}h_2^{(l)}\}$.

Thus, it is proved

**Lemma 2.** Approximations $(V_B^{N_B}, \tilde{p}_B^{N_B}, \tilde{r}_B^{N_B})$ of the space $H^m(\Gamma)$ are convergent and for arbitrary $u \in H^m(\Gamma)$ are valid the estimates

$$\left\| u - \tilde{p}_B^{N_B} \tilde{r}_B^{N_B} u \right\|_{H^t(\Gamma)} \leq C h^{\sigma-t} \| u \|_{H^\sigma(\Gamma)}, \ 0 \leq t \leq \sigma \leq m+1, \ t \leq m, \qquad (12)$$

in which constant $C > 0$ does not depend on $u$.

## 4. Lagrangian approximations

Assign to each element of the grid $S_h$ of domain $S$

$$P_{ij} = [h_1 i, h_1(i+1)] \times [h_2 j, h_2(j+1)], \ i = \overline{0(1)(n-1)}, \ j = \overline{0(1)(k-1)},$$

a smaller rectangular grid $P_{ij}^\varepsilon$ with the steps $\varepsilon_1 = h_1/m$ and $\varepsilon_2 = h_2/m$. Denote $S_{h,\varepsilon} = \bigcup_{i,j} P_{ij}^\varepsilon$ and associate with the set of nodes $S_{h,\varepsilon}$ a system of piecewise polynomial functions

$$\{L_{pt}(\xi)\}_{p=0 \ t=0}^{mn \ \ mk}, \qquad (13)$$

satisfying conditions

$$L_{pt}(\xi_{ls}) = \delta_{pl} \delta_{ts}, \ \mathrm{supp}\{L_{pt}(\xi)\} = \tilde{P}_{pt}, \ \tilde{P}_{pt} = \{\bigcup_{i,j} P_{ij} : \xi_{pt} \in P_{ij}\}, \ \xi_{ls} \in \tilde{P}_{pt}, \qquad (14)$$

where $\delta_{pl}$ is the Kronecker symbol.

Functions (13)-(14) form a system of Lagrangian finite elements of $m$-th degree in $H^m(S)$. Denote by $U_L^{N_1}$ the linear shell of this system, $N_1 = (1+mn)(1+mk)$. It is obvious that the restriction of system (13)-(14) onto an arbitrary rectangle $P_{ij}$ of the grid $S_h$ is a basis in the space of polynomials $\mathrm{P}^m(P_{ij})$ of degree not higher than $m$, defined on $P_{ij}$. Then

$$U_B^N \subset U_L^{N_1}. \qquad (15)$$

Choose the extension operator $p_L^{N_1} : V_L^{N_1} \to U_L^{N_1} \subset H^m(S)$, where $V_L^{N_1} \subset R^{N_1}$, in the form

$$p_L^{N_1} \mathbf{v}_L^{N_1} = \sum_{i=0}^{mn} \sum_{j=0}^{mk} v_{N_1}^{(i,j)} L_{ij}(\xi), \ \mathbf{v}_L^{N_1} = (v_{N_1}^{(0,0)}, v_{N_1}^{(0,1)}, ..., v_{N_1}^{(mn,mk)}). \qquad (16)$$

Then, by virtue of the embedding (15), there exists a restriction operator $r_L^{N_1} : H^m(S) \to V_L^{N_1}$ such that approximations $(V_L^{N_1}, p_L^{N_1}, r_L^{N_1})$ of the space $H^m(S)$ are convergent and valid the estimates

$$\left\| v - p_L^{N_1} r_L^{N_1} v \right\|_{H^t(S)} \leq \tilde{C} h^{\sigma-t} \| v \|_{H^\sigma(S)}, \ 0 \leq t \leq \sigma \leq m+1, \ t \leq m, \qquad (17)$$

in which constant $\tilde{C} > 0$ does not depend on $v$.

Thus, it is proved

**Lemma 3.** There is a restriction operator $r_L^{N_1} : H^m(S) \to V_L^{N_1}$ such that approximations $(V_L^{N_1}, p_L^{N_1}, r_L^{N_1})$ of the space $H^m(S)$ are convergent and valid the estimates (17).

Assume that surface $\Gamma$ satisfy the conditions of p. 3. Construct in each domain $S_l$ the rectangular grid $S_l^h$ with the steps $h_1^{(l)} = a_l/n_l$ and $h_2^{(l)} = b_l/k_l$ and set on each element $P_{ij}^l$ of the grid $S_l^h$ a smaller grid with the steps $\varepsilon_1^{(l)} = h_1^{(l)}/m$ and $\varepsilon_2^{(l)} = h_2^{(l)}/m$, $l = \overline{1,M}$. Define analogously to (13), (14) in each grid domain $S_l^{h,\varepsilon} = \bigcup_{i,j} P_{ij}^{l,\varepsilon}$ the system of Lagrangian finite elements

$$\{L_{ij}^{(l)}(\xi^{(l)})\}_{i=0 \ j=0}^{n_l m \ \ k_l m}, \ \xi^{(l)} \in S_l, \ l = \overline{1,M}.$$

Assign to the family $S_l^{h,\varepsilon}$ the grid $\Gamma_{h,\varepsilon} = \bigcup_{l=1}^{M} \tau_l^{-1}(S_l^{h,\varepsilon})$ on the surface $\Gamma$, where $\tau_l^{-1}(P_{ij}^{l,\varepsilon})$ are the elements of the grid $\Gamma_{h,\varepsilon}$, $l = \overline{1,M}$. Denote by $T_l$ the set of nodes of the grid $S_l^{h,\varepsilon}$, $l = \overline{1,M}$, $T = \bigcup_{l=1}^{M} T_l$. We number all elements of the set $T$ with the cross-cutting index $t = \overline{1,K}$, $K = \sum_{l=1}^{M} K_l$, $K_l = (1 + n_l m)(1 + k_l m)$, and put in correspondence to each node $x_p$ of the grid $\Gamma_{h,\varepsilon}$ the set of elements

$$P_p^* = \{P_{ij}^{l,\varepsilon} \subset \bigcup_{l=1}^{M} S_l^{h,\varepsilon} : x_p \in \tau_l^{-1}(P_{ij}^{l,\varepsilon})\},$$

element

$$\tilde{P}_p = \{\bigcup_{i,j} \tau_l^{-1}(P_{ij}^{l,\varepsilon}), P_{ij}^{l,\varepsilon} \in P_p^*, l = \overline{1,M}\},$$

the set of indexes

$$T_p^* = \{t \in T : \tau_l^{-1}(P_{ij}^{l,\varepsilon}) = x_p, \xi_t^{(l)} \in P_l^{h,\varepsilon} \, l = \overline{1,M}\},$$

and function

$$\tilde{L}_p(x) = \sum_{t \in T_p^*} L_t^l(\tau_l(x)), \, x \in \Gamma_l, \, \mathrm{supp}\{\tilde{L}_p(x)\} = \tilde{P}_p, \, L_t^{(l)}(\xi^{(l)}) \in \{L_{ij}^{(l)}(\xi^{(l)})\}_{i=0}^{n_l m} \, _{j=0}^{k_l m}.$$

Denote by $\tilde{r}_L^{N_L}$ the restriction operator from $H^m(\Gamma)$ into the finite dimensional space $V_L^{N_L}$ and by $\tilde{r}_L^{N_l}$ – its restriction to $H^m(\Gamma_l)$, i.e.

$$\tilde{r}_L^{N_L} = \{\tilde{r}_L^{N_l}\}_{l=1}^{M}, \, \tilde{r}_L^{N_l} u_l(x) = r_L^{K_l} v_l(\xi), \qquad (18)$$

where $r_L^{K_l}$ is the restriction operator from $H^m(S_l)$ into the corresponding finite dimensional space $V_L^{K_l}$, $l = \overline{1,M}$, and $N_L$ is the number of nodes in the grid $\Gamma_{h,\varepsilon}$.

The extraction operator $\tilde{p}_L^{N_L}$ from $V_L^{N_L}$ into the linear shell $U_L^{N_L}$ of the system $\{\tilde{L}_p(x)\}_{p=1}^{N_L}$, $U_L^{N_L} \subset H^m(\Gamma)$, introduce by formula

$$(\tilde{p}_L^{N_L} \mathbf{u}_L^{N_L})(x) = \sum_{i=1}^{N_L} u_N^{(i)} \tilde{L}_i(x), \, \mathbf{u}_L^{N_L} \in V_L^{N_L}. \qquad (19)$$

From Lemma 3 follows that

$$\lim_{N_L \to \infty} \|u - \tilde{p}_L^{N_L} \tilde{r}_L^{N_L} u\|_{H^m(\Gamma)} = \sum_{l=1}^{M} \lim_{N_l \to \infty} \|v_l - p_L^{N_l} r_L^{N_l} v_l\|_{H^m(S_l)} = 0,$$

i.e. approximations $(V_L^{N_L}, \tilde{p}_L^{N_L}, \tilde{r}_L^{N_L})$ of the space $H^m(\Gamma)$ are convergent. Further from estimate (17) we obtain

$$\|u - \tilde{p}_L^{N_L} \tilde{r}_L^{N_L} u\|_{H^t(\Gamma)}^2 = \sum_{l=1}^{M} \|v_l - p_L^{N_l} r_L^{N_l} v_l\|_{H^t(S_l)}^2 \leq \tilde{C}^2 h^{2(\sigma-t)} \sum_{l=1}^{M} \|v_l\|_{H^\sigma(S_l)}^2 = \tilde{C}^2 h^{2(\sigma-t)} \|u\|_{H^\sigma(\Gamma)}^2,$$

$$0 \leq t \leq \sigma \leq m+1, \, t \leq m,$$

where $p_L^{N_l}$ is a similar to (16) extension operator from $V_L^{N_l}$ into $H^m(S_l)$, constant $\tilde{C} > 0$ does not depend on $u$ and $h = \max_{1 \leq l \leq M} \{h_1^{(l)} h_2^{(l)}\}$.

Thus, it is proved

**Lemma 4.** There is a restriction operator $\tilde{r}_L^{N_L} : H^m(\Gamma) \to V_L^{N_L}$ such that approximations $(V_L^{N_L}, \tilde{p}_L^{N_L}, \tilde{r}_L^{N_L})$ of the space $H^m(\Gamma)$ are convergent and valid the estimates

$$\left\| u - \tilde{p}_L^{N_L} \tilde{r}_L^{N_L} u \right\|_{H^t(\Gamma)} \leq \tilde{C} h^{\sigma - t} \|u\|_{H^\sigma(\Gamma)}, \quad 0 \leq t \leq \sigma \leq m+1, \ t \leq m, \quad (20)$$

in which constant $\tilde{C} > 0$ does not depend on $u$.

## 5. Galerkin method

Let us denote $G' = R^3 \setminus \overline{G}$ and introduce in $G$ and $G'$ the Sobolev spaces [30]

$$H^m(G) = \{v \in L_2(G) : \partial^\alpha v \in L_2(G), |\alpha| \leq m\},$$

$$W^m(G') = \{v \in D'(G') : (1 + r^2)^{(|\alpha|-1)/2} \partial^\alpha v \in L_2(G'), |\alpha| \leq m\}.$$

where $m \geq 0$, and $r = (\sum_{i=1}^{3} x_i^2)^{1/2}$, $x = (x_1, x_2, x_3) \in R^3$.

Consider the next boundary value problem: to find function

$$v \in H^{m+1}_{\Gamma, \Delta = 0} = \{v \in H^{m+1}(G) \cup W^{m+1}(G') : v|_{\Gamma_{int}} = v|_{\Gamma_{ext}}, \Delta v(x) = 0, x \in G, G'\} \quad (21)$$

satisfying condition

$$v|_\Gamma = f, \ f \in H^{m+1/2}(\Gamma). \quad (22)$$

In [9] was proved the next

**Theorem 2.** Problem (21)-(22) has one and only one solution.

We will search a solution of the problem (21) - (22) in the form of simple layer potential

$$v(x) = \frac{1}{4\pi} \int_\Gamma \frac{u(y)}{|x-y|} d\Gamma_y, \ x \in G, G'.$$

The unknown potential density is determined from the equation

$$(Au)(x) \equiv \frac{1}{4\pi} \int_\Gamma \frac{u(y)}{|x-y|} d\Gamma_y = f(x), \ x \in \Gamma. \quad (23)$$

The next result is in order [9].

**Theorem 3.** Operator A is an isomorphism of $H^s(\Gamma)$ onto $H^{s+1}(\Gamma)$.

From the last statement and the Banach theorem follows the validity of inequalities

$$\alpha_s \|u\|_{H^s(\Gamma)} \leq \|Au\|_{H^{s+1}(\Gamma)} \leq \beta_s \|u\|_{H^s(\Gamma)}, \quad (24)$$

in which constants $\alpha_s$ and $\beta_s$, $0 < \alpha_s \leq \beta_s$, does not depend on $u \in H^s(\Gamma)$.

Suppose that for approximation of unknown potential density $u \in H^m(\Gamma)$ uses the system of B-splines of the form (9), and $U_{N_B}$ is its linear shell. We choose the operators $\tilde{r}_{N_B} : H^m(\Gamma) \to V_{N_B}$ and $\tilde{p}_{N_B} : V_{N_B} \to U_{N_B}$ in the form (10) and (11) respectively and determine the restriction operator $s_{N_B} : H^{m+1}(\Gamma) \to \Phi_{N_B}$ in the form $s_{N_B} = \tilde{r}_{N_B}$. In this case, the system

$$\mathbf{A}_{N_B}^G \mathbf{u}_{N_B} = \mathbf{f}_{N_B}, \ \mathbf{A}_{N_B}^G = \tilde{r}_{N_B} A \tilde{p}_{N_B}, \ \mathbf{f}_{N_B} = r_{N_B} f,$$

implements Galerkin method of solving the equation (23). From Lax-Milgram lemma [31] follows that matrix $\mathbf{A}_{N_B}^G$ is nondegenerate and, accordingly, the definition of operator $q_{N_B}$ in the form $q_{N_B} \mathbf{f}_{N_B} = A \tilde{p}_{N_B} \mathbf{u}_{N_B}$ is correct. Taking into account the left side of inequalities (24), the bijectivity of mapping $\tilde{p}_{N_B} : V_{N_B} \to U_{N_B}$, the expressions for the norms in the spaces $V_{N_B}$ and $\Phi_{N_B}$ in the case $U = H^m(\Gamma)$, $F = H^{m+1}(\Gamma)$, and equality $Q_{N_B} A P_{N_B} u = A P_{N_B} u$, we obtain the following inequalities

$$\alpha_m \|\mathbf{u}_{N_B}\|_{V_{N_B}} \leq \|\mathbf{A}_{N_B}^G \mathbf{u}_{N_B}\|_{\Phi_{N_B}} \quad (25)$$

for arbitrary $\mathbf{u}_{N_B} \in V_{N_B}$ and $\alpha_m$ does not depend on $\mathbf{u}_{N_B}$.

Then from the inequalities (24) and (25), Lemma 2 and Theorem 1 we obtain the validity of following statement.

Theorem 4. For arbitrary $f \in H^{m+1}(\Gamma)$, $m=0,1,\ldots$, the approximate solution $u_{N_B}^B$ of equation (23) obtained by the Galerkin method under approximation of unknown potential density by the system of functions constructed on the basis of B-splines of $m$-th degree converges to its exact solution, and there is an estimate

$$\|u - u_{N_B}^B\|_{H^t(\Gamma)} \leq \frac{C(1+\beta_t/\alpha_t)}{\alpha_\sigma} h^{\sigma-t} \|f\|_{H^{\sigma+1}(\Gamma)}, \quad 0 \leq t \leq \sigma \leq m+1, \ t \leq m, \quad (26)$$

where h is the maximum area of the grid element on $\Gamma$.

Similarly, from the inequalities (24) and (25), Lemma 4 and Theorem 1, we obtain the validity of following statement.

Theorem 5. For arbitrary $f \in H^{m+1}(\Gamma)$, $m=0,1,\ldots$, the approximate solution $u_{N_L}^L$ of equation (23) obtained by the Galerkin method under approximation of unknown potential density by the system of functions constructed on the basis of Lagrangian finite elements of $m$-th degree converges to its exact solution, and there is an estimate

$$\|u - u_{N_L}^L\|_{H^t(\Gamma)} \leq \frac{\tilde{C}(1+\beta_t/\alpha_t)}{\alpha_\sigma} h^{\sigma-t} \|f\|_{H^{\sigma+1}(\Gamma)}, \quad 0 \leq t \leq \sigma \leq m+1, \ t \leq m, \quad (27)$$

where h is the maximum area of the grid element on $\Gamma$.

## 6. Collocation method

To simplify the presentation, we assume that for approximation of unknown potential density $u \in H^m(\Gamma)$, $m \geq 0$, of equation (23) a system of linearly independent functions $\{\varphi_i\}_{i=1}^\infty$ is chosen, $U_N$ is a linear shell of the system $\{\varphi_i\}_{i=1}^N$, $r_N : H^m(\Gamma) \to V_N$, $p_N : V_N \to U_N$ are the similar to described in p. 2 restriction and extraction operators. Denote by $X_N$ the set of pairwise different points belonging to the surface $\Gamma$

$$X_N = \{x_j\}_{j=1}^N, \ x_j \in \Gamma, \ j = \overline{1,N},$$

and introduce in $H^{m+1}(\Gamma)$ restriction operator $s_N : H^{m+1}(\Gamma) \to \Phi_N$ by formula

$$(s_N f)_j = f(\tilde{y}_j) \quad (28)$$

in which

$$\tilde{y}_j \in \{\tilde{y} \in \delta(y_j) : |f(\tilde{y})| = \min_{y \in \delta(y_j)} |f(y)|, \ y_j \in X_N\}, \ \delta(y_j) = \{y \in \Gamma : |y - y_j| < \delta\}, \quad (29)$$

in particular $\rho(y^*, \delta(y_j)) > 0$ for arbitrary $y^* \in X_N$, $y^* \neq y_j$, $j = \overline{1,N}$.

If $f \in C(\Gamma)$, then operator $s_N$ can be defined as usual

$$(s_N f)_j = f(y_j), \ y_j \in X_N, \quad (30)$$

i.e. $\tilde{y}_j = y_j$, $j = \overline{1,N}$. It is easy to see that, with this choice of operator $s_N$, a system of linear algebraic equations

$$\mathbf{A}_N^c \mathbf{u}_N = s_N f, \ \mathbf{A}_N^c = s_N A p_N, \ \mathbf{u}_N \in V_N, \quad (31)$$

implements the collocation method of solving the equation (23). The set $X_N$ is called a set of collocation points.

Denote $Y_N = \{\tilde{y}_j\}_{j=1}^N$ and consider the system of functions

$$r_j(x) = \frac{1}{|x - \tilde{y}_j|}, \quad \tilde{y}_j \in Y_N, \quad j = \overline{1, N}.$$

From the choice of the set $X_N$ and conditions (29) follow that the functions of system $\{r_j(x)\}_{j=1}^N$ are linearly independent [32].

Define in $L^\infty(\Gamma)$ the family of linear continuous functionals

$$l_j(\varphi) = \int_\Gamma \varphi(x) r_j(x) d\Gamma_x, \quad \varphi \in L^\infty(\Gamma), \quad j = \overline{1, N}.$$

Denote by $Ker(l_j)$ the zero subspace of functional $l_j$ in $L^\infty(\Gamma)$

$$Ker(l_j) = \{\varphi \in L^\infty(\Gamma) : l_j(\varphi) = 0\}$$

and suppose that $K_N = \bigcap_{j=1}^N Ker(l_j)$. The degeneracy of matrix $\mathbf{A}_N^c$ is equivalent to the linear dependence of its rows or columns, that is, the existence of such sets $\boldsymbol{\alpha}_N = \{\alpha_i\}_{i=1}^N \in R^N$ or $\boldsymbol{\beta}_N = \{\beta_j\}_{j=1}^N \in R^N$, $\sum_{i=1}^N \alpha_i^2 > 0$, $\sum_{j=1}^N \beta_j^2 > 0$, that

$$\int_\Gamma (\sum_{i=1}^N \alpha_i \varphi_i(x)) r_j(x) d\Gamma_x = 0, \quad j = \overline{1, N}, \tag{32}$$

or

$$\int_\Gamma \varphi_i(x) (\sum_{j=1}^N \beta_j r_j(x)) d\Gamma_x = 0, \quad i = \overline{1, N}. \tag{33}$$

Implementation of equations (32), (33) is only possible if $K_N \cap U_N \neq 0$. From this follows sufficient conditions for the invertibility of matrix $\mathbf{A}_N^c$, which we formulate in the next statement.

Lemma 5. Let us the system of linearly independent functions $\{\varphi_i\}_{i=1}^N$ is chosen for the approximate solution of equality (23) and determined the set of collocation points $X_N$ (and, consequently, the set $K_N$ is defined). Then, if

$$K_N \cap U_N = 0, \tag{34}$$

then the matrix $\mathbf{A}_N^c$ of the system of collocation equations (31) is non-degenerate for arbitrary $N$.

A similar result is obtained if the restriction operator $s_N$ is chosen in the form

$$(s_N f)_j = \frac{1}{mes\, \delta(y_j)} \int_{\delta(y_j)} f(y) d\Gamma_y \tag{35}$$

and

$$r_j(x) = \frac{1}{mes\, \delta(y_j)} \int_{\delta(y_j)} \frac{d\Gamma_y}{|x - y|}, \quad j = \overline{1, N}.$$

It is obvious that under conditions of Lemma 5 the operator $\mathbf{A}_N^c$, where $s_N$ is defined according to (27)-(28) or (33), or in the case of $f \in C(\Gamma)$ according to (30), is stable in sense (4).

Consider a discrete analog of condition (34). Let us the quadrature formula

$$\int_\Gamma \varphi(x) r_i(x) d\Gamma_x \approx \sum_{j=1}^{N} A_j \varphi(x_j) r_i(x_j), \ x_j \in \Gamma, \ x_j \neq x_i, \ \text{if} \ j \neq i, \qquad (36)$$

is used to calculate the integrals

$$\int_\Gamma \varphi(x) r_i(x) d\Gamma_x, \ \varphi(x) \in U_N, \ i = \overline{1,N},$$

which is exact for integrals

$$\int_\Gamma \varphi(x) \psi(x) d\Gamma_x, \ \varphi(x), \psi(x) \in U_N.$$

Consider the system of functions

$$\psi_i(x) = \sum_{k=1}^{N} \alpha_k^{(i)} \varphi_k(x), \qquad (37)$$

the coefficients $\alpha_k^{(i)}$, $k, i = \overline{1,N}$, of which we define from $N$ systems of linear algebraic equations

$$\sum_{k=1}^{N} \alpha_k^{(i)} \varphi_k(x_j) = r_i(x_j), \ i, j = \overline{1,N}. \qquad (38)$$

Define the conditions under which the functions $\psi_i(x)$, $i = \overline{1,N}$, are linearly independent. From (37) we obtain that

$$\sum_{i=1}^{N} c_i \psi_i(x) = \sum_{k=1}^{N} (\sum_{i=1}^{N} c_i \alpha_k^{(i)}) \varphi_k(x) = 0$$

if and only if

$$\sum_{i=1}^{N} c_i \alpha_k^{(i)} = 0, \ k = \overline{1,N}. \qquad (39)$$

Let us the set of colocation points $X_N = \{y_j\}_{j=1}^{N} \subset \Gamma$ is chosen in such a way that

$$0 < |x_i - y_i| < \varepsilon, \ d < |x_i - y_j|, \ i, j = \overline{1,N}, \ j \neq i, \ 0 < \varepsilon < \frac{d}{N-1},$$

where $\{x_j\}_{j=1}^{N}$ are the nodes of quadrature formula (36). Then

$$r_i(x_i) > \sum_{i=1, i \neq j}^{N} r_i(x_j),$$

matrix $\mathbf{R}_N = \{r_i(x_j)\}_{i,j=1}^{N}$ due to Hadamard condition is nondegenerate and from (38) we obtain that vectors $\boldsymbol{\alpha}_k = \{\alpha_k^{(j)}\}_{j=1}^{N}$, $k = \overline{1,N}$, are linearly independent. Hence, equality (39) holds if and only if $c_i = 0$, $i = \overline{1,N}$, i.e. the functions of system $\{\psi_i(x)\}_{i=1}^{N}$ are linearly independent.

Now, if the quadrature formula of form (36) is used to calculate the integrals in coefficients of matrix $\mathbf{A}_N^c$, instead of the system of collocation equations (31), we actually solve a system with matrix

$$\widetilde{\mathbf{A}}_N^c = \{\int_\Gamma \varphi_i(x) \psi_j(x) d\Gamma_x\}_{i,j=1}^{N},$$

where functions $\psi_i(x)$, $i = \overline{1,N}$, are defined by formulas (37) and (38). The last matrix can be degenerate if and only if there exists a nonzero element $\varphi(x) = \sum_{i=1}^{N} a_i \varphi_i(x) \in U_N$, orthogonal to all $\psi_i(x)$, $i = \overline{1,N}$, which is impossible, since the system $\{\psi_i(x)\}_{i=1}^{N}$ forms a basis in the space $U_N$.

Let us the system of B-splines of the form (9) is used to approximate the unknown potential density $u \in H^m(\Gamma)$ and $U_{N_B}$ is its linear shell. We choose the operators $\tilde{r}_{N_B} : H^m(\Gamma) \to V_{N_B}$ and $\tilde{p}_{N_B} : V_{N_B} \to U_{N_B}$ in the form (18) and (19) respectively and determine the restriction operator . $s_{N_B} : H^{m+1}(\Gamma) \to \Phi_{N_B}$ in the form (28), (29). In this case, the system

$$\mathbf{A}^c_{N_B} \mathbf{u}_{N_B} = \mathbf{f}_{N_B}, \quad \mathbf{A}^c_{N_B} = \tilde{r}_{N_B} A \tilde{p}_{N_B}, \quad \mathbf{f}_{N_B} = r_{N_B} f,$$

implements the collocation method for solution of equation (23). From Lax-Milgram lemma [31] follows that under conditions (34) matrix $\mathbf{A}^c_{N_B}$ is non-degenerate and, accordingly, the definition of operator $q_{N_B}$ in the form $q_{N_B} \mathbf{f}_{N_B} = A \tilde{p}_{N_B} \mathbf{u}_{N_B}$ is correct. Given the left side of inequalities (24), the biectivity of mapping $\tilde{p}_{N_B} : V_{N_B} \to U_{N_B}$, the expressions for norms in the spaces $V_{N_B}$ and $\Phi_{N_B}$ in the case $U = H^m(\Gamma)$, $F = H^{m+1}(\Gamma)$ and equality $Q_{N_B} A P_{N_B} u = A P_{N_B} u$, we obtain the validity of inequalities (25) for arbitrary $\mathbf{u}_{N_B} \in V_{N_B}$, in which $\alpha_m$ does not depend on $\mathbf{u}_{N_B}$.

Then from the inequalities (24) and (25), Lemmas 2, 5, and Theorem 1 we obtain the validity of following statement.

**Theorem 6.** For arbitrary $f \in H^{m+1}(\Gamma)$, $m=0,1,\ldots$, the approximate solution $u^B_{N_B}$ of equation (23) obtained by collocation method under approximation of unknown potential density by a system of functions constructed on the basis of B-splines of $m$-th degree and the choice of collocation points that satisfies the condition (34) converges to its exact solution, and there is an estimate

$$\left\| u - u^B_{N_B} \right\|_{H^t(\Gamma)} \leq \frac{C(1 + \beta_t / \alpha_t)}{\alpha_\sigma} h^{\sigma - t} \| f \|_{H^{\sigma+1}(\Gamma)}, \quad 0 \leq t \leq \sigma \leq m+1, \ t \leq m, \qquad (40)$$

where h is the maximum area of the grid element on $\Gamma$.

Similarly, from the inequalities (24) and (25), Lemmas 4, 5, and Theorem 1 we obtain the validity of following statement.

**Theorem 7.** For arbitrary $f \in H^{m+1}(\Gamma)$, $m=0,1,\ldots$, the approximate solution $u^L_{N_L}$ of equation (23) obtained by collocation method under approximation of unknown potential density by a system of functions constructed on the basis of Lagrangian finite elements of $m$-th degree and the choice of collocation points that satisfies the condition (34) converges to its exact solution, and there is an estimate

$$\left\| u - u^L_{N_L} \right\|_{H^t(\Gamma)} \leq \frac{\tilde{C}(1 + \beta_t / \alpha_t)}{\alpha_\sigma} h^{\sigma - t} \| f \|_{H^{\sigma+1}(\Gamma)}, \quad 0 \leq t \leq \sigma \leq m+1, \ t \leq m, \qquad (41)$$

where h is the maximum area of the grid element on $\Gamma$.

## 7. Error estimation of approximate solution of the Dirichlet problem for the Laplace equation

Denote by $u_N(x)$ the approximate solution of equation (23), obtained by means of considered above Galerkin or collocation methods, $N = N_B$ in the case of approximation by B-splines and $N = N_L$ in the case of Lagrangian approximations. Denote

$$v_N(x) = \frac{1}{4\pi} \int_\Gamma \frac{u_N(y)}{|x - y|} d\Gamma_y, \quad x \in G, G',$$

and estimate the modulus of value

$$\frac{\partial^\alpha}{\partial x^\alpha}(v(x) - v_N(x)) = \frac{1}{4\pi} \int_\Gamma (u(y) - u_N(y)) \frac{\partial^\alpha}{\partial x^\alpha} \frac{1}{|x - y|} d\Gamma_y, \quad x \in G, G', \alpha = 0,1,\ldots$$

Let us

$$x \in R^3 \setminus \{\tilde{x} \in R^3 : |\tilde{x} - y| < \delta, y \in \Gamma\}. \qquad (42)$$

Using Holder inequality, we obtain

$$\left|\frac{\partial^\alpha}{\partial x^\alpha}(v(x)-v_N(x))\right| \leq \|u-u_N\|_{L_2(\Gamma)} \sqrt{\int_\Gamma \frac{\partial^\alpha}{\partial x^\alpha}\frac{1}{|x-y|}d\Gamma_y}, \quad x \in G, G',$$

or, taking into account (42),

$$\left|\frac{\partial^\alpha}{\partial x^\alpha}(v(x)-v_N(x))\right| \leq \frac{mes\Gamma}{\delta^{\alpha+1}}\|u-u_N\|_{L_2(\Gamma)}, \quad x \in G, G', \quad \alpha = 0,1,\ldots \qquad (43)$$

Then from inequalities (24), (43) and Theorems 4-7 follow the validity of the next statement.

Theorem 8. For arbitrary $f \in H^{m+1}(\Gamma)$, $m=0,1,\ldots$, an approximate solution of the problem (21), (22) obtained by Galerkin or collocation methods under approximation of unknown potential density by systems of functions constructed on the basis of B-splines or Lagrangian finite elements of the $m$-th degree, converges to its exact solution, and there is an estimate

$$\left|\frac{\partial^\alpha}{\partial x^\alpha}(v(x)-v_N(x))\right| \leq \frac{C^*(1+\beta_0/\alpha_0)h^m}{\alpha_m \delta^{\alpha+1}}\|f\|_{H^{m+1}(\Gamma)}, \quad x \in G, G', \quad \alpha = 0,1,\ldots$$

**Conclusions**

The paper describes the conditions and evaluations of convergence of Galerkin and collocation methods for solution of Fredholm integral equation of the first kind for the simple layer potential in case of closed boundary surface in a three-dimensional space. Approximation of potential density was performed using B-splines and Lagrangian finite elements of various orders on rectangular grids constructed in the desired function definition domain. Estimations were obtained for the error of approximate solution of Dirichlet problem for Laplace equation that is equivalent to the integral equation for the simple layer potential. The approach proposed can be used to define convergence of other projection methods (the smallest squares, smallest mismatch etc.) for solving potential theory integral equations that are equivalent to the boundary value problems for equations of mathematical physics and other types of finite elements of various orders, constructed on both rectangular and triangular grids in desired potential density definition domain.